\documentclass[10pt]{amsart}
\usepackage{graphicx}
\usepackage[centertags]{amsmath}
\usepackage{amsfonts}
\usepackage{amssymb}
\usepackage{amsthm}
\usepackage{newlfont}
\newtheorem{theorem}{Theorem}[section]

\newtheorem{lemma}[theorem]{Lemma}

\theoremstyle{definition}

\theoremstyle{remark}

\newcommand{\id}{\textrm{id}}

\title{The value of the work done by an isotropic vector force field along an isotropic curve}

\author{Dimitar Razpopov$^{1}$ and Georgi Dzhelepov$^{2}$}

\address{$^{1}$Department of Mathematics and Informatics, Agricultural University of Plovdiv,
12 Mendeleev Blvd., Bulgaria 4000}

\email{razpopov@au-plovdiv.bg}

\address{$^{2}$Department of Mathematics and Informatics, Agricultural University of Plovdiv,
12 Mendeleev Blvd., Bulgaria 4000}

\email{dzhelepov@au-plovdiv.bg}

\begin{document}

\begin{abstract}
In the present paper we consider a 3-dimensional differentiable manifold $M$ equipped with a Riemannian metric $g$ and an endomorphism $Q$, whose third power is the identity and $Q$ acts as an isometry on $g$. Both structures $g$ and $Q$ determine an associated metric $f$ on $(M, g, Q)$. The metric $f$ is necessary indefinite and it defines isotropic vectors in the tangent space $T_{p}M$ at an arbitrary point $p$ on $M$.

The physical forces are represented by vector fields. We investigate physical forces whose vectors are in $T_{p}M$ on $(M, g, Q)$. Moreover, these vectors are isotropic and they act along isotropic curves. We study the physical work done by such forces.
\end{abstract}



\maketitle

\section{Introduction}
The physical work and the physical force on differentiable manifolds have a great application in physics.  Vector fields are often used to model a force, such as the magnetic or gravitational force, as it changes from one point to another point. As a particle moves through a force field along a curve $c$, the work done by the force is the product of force and displacement. There are some papers concerning physical results on light-like (degenerate) objects of differentiable manifolds (\cite{duggal}, \cite{Korpinar} and \cite{robert}).

The object of the present paper is a 3-dimensional differentiable manifold $M$ equipped with a Riemannian metric $g$ and a tensor $Q$ of type $(1, 1)$, whose third power is the identity and $Q$ acts as an isometry on $g$. Such a manifold $(M, g, Q)$ is defined in \cite{dzhelepov2} and studied in  \cite{AE-dok}, \cite{fil-dok}, \cite{dok-raz-dzhe} and \cite{dzhelepov1}. Also, we consider an associated metric $f$,  which is introduced in \cite{dzhelepov1}. The metric $f$ is necessary indefinite and it determines space-like vectors, isotropic vectors and time-like vectors in the tangent space $T_{p}M$ at an arbitrary point $p$ on $M$.

 We investigate physical forces whose vectors are in $T_{p}M$ on $(M, g, Q)$. Moreover, these vectors are isotropic with respect to $f$ and they act along isotropic curves. We study the physical work done by such forces.

\section{Preliminaries}
Let $M$ be a $3$-dimensional Riemannian manifold equipped with an endomorphism $Q$ in the tangent space $T_{p}M$, $p \in M$. Let the local coordinates of $Q$ with respect to some coordinate system form the circulant matrix:
\begin{equation*}%\label{f2}
    (Q_{i}^{j})=\begin{pmatrix}
      0 & 1 & 0 \\
      0 & 0 & 1 \\
      1 & 0 & 0 \\
    \end{pmatrix}.
\end{equation*}
Then $Q$ has the property
\begin{equation}\label{q3}
    Q^{3}=\id.
\end{equation}
Let $g$ be a positive definite metric on $M$, which satisfies the equality
\begin{equation}\label{ff1}
 \quad g(Qr, Qi)=g(r,i).
 \end{equation}
 In \eqref{ff1} and further $r, i, w$ will stand for arbitrary vectors in $T_{p}M$.

 Such a manifold $(M, g, Q)$ is introduced in \cite{dzhelepov2}.

 It is well-known that the norm of every vector $i$ is given by
 $\|i\|=\sqrt{g(i,i)}.$
 Then, having in mind (\ref{ff1}), for the angle $\varphi = \angle(i,Qi)$ we have
 \begin{equation*}%\label{ff3}
 \cos \varphi=\frac{g(i,Qi)}{g(i,i)}.
 \end{equation*}
 In \cite{dzhelepov2}, for  $(M, g, Q)$, it is verified that the angle $\varphi$ is in $[0,\frac{2\pi}{3}]$. If $\varphi \in (0, \frac{2\pi}{3})$, then the vector $i$ form a basis $\{i,Qi, Q^{2}i\}$, which is called a $Q$-basis of $T_{p}M$.

 The associated metric $f$ on $(M,g,Q)$, determined by
\begin{equation}\label{f3}
   f(r,i)=g(r,Qi)+g(Qr,i).
\end{equation}
is necessary indefinite \cite{dzhelepov1}.

A vector $r$ in $T_{p}M$ is isotropic with respect to $f$ if
\begin{equation}\label{f4}
f(r,r)=0.
\end{equation}

In every $T_{p}M$, for $(M, g, Q)$, there exists an orthonormal $Q$-basis $\{i, Qi, Q^{2}i\}$ (\cite{dzhelepov2}).
From \eqref{q3}, (\ref{f3}) and (\ref{f4}) we state the following
\begin{lemma}
Let $\{i, Qi, Q^{2}i\}$ be an orthonormal $Q$-basis of $T_{p}M$. If $r=ui+vQi+qQ^{2}i$ is an isotropic vector, then its coordinates satisfy
\begin{equation}\label{f5}
uv+vq+qu=0.
\end{equation}
\end{lemma}

An isotropic (null) curves $c: r=r(t)$ are those whose tangent vectors are everywhere isotropic, i.e.,
\begin{equation}\label{is-curve}
f(dr, dr)=0.
\end{equation}

The physical forces are represented by vector fields. We investigate physical forces whose vectors are in $T_{p}M$ on $(M, g, Q)$. Moreover, these vectors are isotropic and they act along isotropic curves. We study the physical work done by such forces.
%%% ----------------------------------------------------------------

\section{The work in $T_{p}M$}
We consider an orthonormal $Q$-basis $\{i, Qi, Q^{2}i\}$ in $T_{p}M$ on $(M, g, Q)$.

Let $p_{xyz}$ be a coordinate system such that the vectors $i$, $Qi$ and $Q^{2}i$ are on the axes $p_{x}$, $p_{y}$ and $p_{z}$, respectively. So $p_{xyz}$ is an orthonormal coordinate system.

The curve $c$ is determined by
\begin{equation}\label{ff4}
c:r(t)=x(t)i+y(t)Qi+z(t)Q^{2}i,
\end{equation}
where  $t\in [\alpha, \beta]\subset \mathbb{R}$.

Let $c$ be an isotropic smooth curve. Thus equalities \eqref{q3}, \eqref{f3},  \eqref{is-curve} and \eqref{ff4} imply
\begin{equation}\label{f7}
dxdy+dydz+dxdz=0.
\end{equation}

We determine a vector force field
\begin{equation}\label{force}
 F(x,y,z)=P(x,y,z)i+R(x,y,z)Qi+S(x,y,z)Q^{2}i,
 \end{equation}
where $P=P(x, y, z)$, $R=R(x, y, z)$, $S=S(x, y, z)$ are smooth functions.

Let the vector field $F$ be isotropic. Hence following \eqref{f5} we get
\begin{equation}\label{f6}
PR+RS+SP=0.
\end{equation}

Work $A$ done by a force $F$, with respect to $f$, moving along a curve $c$ is given by
\begin{equation}\label{f8}
A=\int_{c}f(F,dr),
\end{equation}
where
\begin{equation}\label{f9}
dr=dxi+dyQi+dzQ^{2}i.
\end{equation}

Case (A) Let $F$ and $c$ are both isotropic and they are on the same direction. Since $c$ is a trajectory of $F$ we have that the vectors $F$ and $dr$ are collinear. Therefore their coordinates satisfy
\begin{equation}\label{f10}
\frac{dx}{P}=\frac{dy}{R}=\frac{dz}{S}=\frac{1}{k},
\end{equation}
where  $k\neq0$ is a function. From \eqref{force} and \eqref{f10} it follows $F=kdr$. Then, having in mind \eqref{is-curve} and \eqref{f8}, we get $dA =f(kdr,dr)= kf(dr,dr)=0$, i.e., $A=0$.

Case (B) Now, we consider the case when $F$ and $c$ are both isotropic but they are on different directions.

From \eqref{f3}, \eqref{f8} and \eqref{f9} it follows
\begin{equation}\label{work}
A=\int_{c}\big[P(dy+dz)+R(dx+dz)+S(dx+dy)],
\end{equation}
and hence
\begin{equation}\label{work2}
A=\int_{\alpha}^{\beta}\big[P(y'+z')+R(x'+z')+S(x'+y')]dt.
\end{equation}

\begin{itemize}
  \item Let $dx+dy=0$. From (\ref{f7}) we have $dx=dy=0$ and $dz\neq 0$. Then $dr=dzQ^{2}i$. Therefore, using \eqref{work2},
  we get
  \begin{equation}\label{f10*}
  A=\int^{\beta}_{\alpha}\big(P(k_{1},k_{2},t)+R(k_{1},k_{2},t)\big)dt,
  \end{equation}
  where $k_{1}$ and $k_{2}$ are specific constants.
  \item Let $P+R=0$. From \eqref{f6} we have $P=R=0$ and hence $S\neq 0$. In this case equalities \eqref{f8} and \eqref{work2} imply
  \begin{equation}\label{f11}
  A = \int^{\beta}_{\alpha}S(x(t),y(t),z(t))[x'(t)+y'(t)]dt.
  \end{equation}
  \item Let $dx+dy\neq 0$ and $P+R\neq0$. With the help of \eqref{f7} and \eqref{f6} we get
  \begin{equation}\label{f12}
  dz=-\frac{dxdy}{dx+dy}, \quad S=-\frac{PR}{P+R}.
  \end{equation}
We use \eqref{work}, \eqref{work2} and \eqref{f12} and obtain that the work $A$ is determined by
  \begin{equation}\label{f13}
  A = \int^{\beta}_{\alpha}\frac{(Py'-Rx')^{2}}{(P+R)(x'+y')}dt.
  \end{equation}

\end{itemize}

From Case (A) and Case (B) we state the following
\begin{theorem}
Let $f$ be the associated metric on $(M,g,Q)$. Let $p_{xyz}$ be a coordinate system such that the vectors $i$, $Qi$ and $Q^{2}i$ of the orthonormal $Q$-basis in $T_{p}M$ are on the axes $p_{x}$, $p_{y}$ and $p_{z}$, respectively. Let $F$ be an isotropic vector force field moving along an isotropic curve $c$. Let $A$ be the work done by $F$. Then
\begin{itemize}
\item[(i)] $A$ is zero if $F$ and $c$ are on the same direction;

\item[(ii)]   $A$ is \eqref{f10*} if $F$ and $c$ are on different directions and $dx+dy=0$;

\item[(iii)] $A$ is  \eqref{f11} if $F$ and $c$ are on different directions and $P+Q=0$;

\item[(iv)] $A$ is  \eqref{f13} if $F$ and $c$ are on different directions and $dx+dy\neq 0$ and $P+R\neq0$.
\end{itemize}
\end{theorem}
\section{Work in a $2$-plane}
Now we consider an arbitrary $2$-plane $\alpha=\{i,Qi\}$ in $T_{p}M$.
We suppose that the angle $\varphi=\angle(i,Qi)$ belongs to the interval $(0, \frac{2\pi}{3}]$.
On $\alpha$ we construct a coordinate system $p_{xy}$ such that $i$ is on the axis $p_{x}$ and $j$ is on the axis $p_{y}$, where
\begin{equation} \label{defw}
j=\frac{1}{\sin \varphi}(-\cos \varphi i+Qi).
\end{equation}
We assume that $\|i\|=1$ and then $p_{xy}$ is an orthonormal coordinate system.

In \cite{sf} it is proved the following
\begin{theorem}
Let $f$ be the associated metric on $(M,g,Q)$ and let $\alpha=\{i, Qi\}$ be an arbitrary $2$-plane in $T_{p}M$. Let the vector $j$ be defined by \eqref{defw} and $p_{xy}$ be a coordinate system such that $i\in p_{x}$, $j\in p_{y}$.  Then the equation of the circle $ c:\ f(w, w)=a^{2}$ in $\alpha$ is given by
\begin{equation}\label{18}
 (\cos\varphi)x^{2}+\frac{(1-\cos\varphi)(1+2\cos\varphi)}{\sin\varphi}xy-\frac{\cos^{2}\varphi}{1+\cos\varphi}y^{2}=\frac{a^{2}}{2}\ ,
\end{equation}
where $\varphi\in(0, \frac{2\pi}{3}]$.
\end{theorem}

Let $w=ui+vj$ be an isotropic vector, i.e., $f(w,w)=0$.
Therefore, with the help of \eqref{18}, we obtain
\begin{equation}\label{ff5}
\cos^{2}\varphi\big(\frac{y}{x}\big)^{2}-\sin\varphi(1+2\cos\varphi)\frac{y}{x}-(1+\cos\varphi)\cos\varphi=0.
\end{equation}
The discriminant of (\ref{ff5}) is
\begin{equation*}%\label{f16}
D = (1+\cos \varphi)(1+3\cos \varphi).
\end{equation*}
Then we get the following cases:

 \textbf{Case (A)} If $\varphi \in(\arccos (-\frac{1}{3}), \frac{2\pi}{3})$, then $D<0$. There is no isotropic directions in $T_{p}M$.

 \textbf{Case (B)} If $\varphi=\arccos(-\frac{1}{3})$, then $D=0$. We have one isotropic straight line $c:\ y=\sqrt{2}x$.
 Then the force $F$ and the curve $c$ both are on one isotropic direction and the work $A$ of the force $F$ along $c$ is zero.

\textbf{ Case (C)} If $\varphi\in\big(0, \frac{\pi}{2}\big)\bigcup\big(\frac{\pi}{2}, \arccos(-\frac{1}{3})\big)$, then $D>0$. We have two isotropic directions which generate two straight lines:
  \begin{equation*}%\label{f17}
 c_{1}: y = k_{1}x, \quad c_{2}: y=k_{2}x, \quad x\in[\alpha,\beta],
  \end{equation*}
where $k_{1}$ and $k_{2}$ are solutions of the equation (\ref{ff5}) for $\frac{y}{x}$.
\begin{itemize}
  \item If $F$ is on $c_{1}$, then the work of $F$ along $c_{1}$ is zero. Similarly, if $F$ is on $c_{2}$, then the work of $F$ along $c_{2}$ is zero.
  \item We suppose that $F$ is on $c_{2}$ but $F$ acts on $c_{1}$. Then
  \begin{equation}\label{ff6}
  F(x,y)=P(x,y)(i+k_{2}j), \quad dr=dt(i+k_{1}j).
  \end{equation}
 Bearing in mind \eqref{ff1} and \eqref{defw} we calculate
 \begin{equation}\label{ff7}
 \begin{array}{ll}
 g(i,Qi)=g(Qi,i)=\cos\varphi, & g(i,Qj)=g(Qj,i)=\frac{\cos\varphi-\cos^{2}\varphi}{\sin\varphi},\\ g(j,Qi)=g(Qi,j)=\sin\varphi, & g(j,Qj)=g(Qj,j)=-\frac{\cos^{2}\varphi}{1+\cos\varphi}.
 \end{array}
 \end{equation}
On the other hand the solutions $k_{1}$ and $k_{2}$ of \eqref{ff5} satisfy equalities
 \begin{equation}\label{ff9}
 k_{1}+k_{2}=\frac{\sin\varphi(1+2\cos\varphi)}{\cos^{2}\varphi}, \quad k_{1}k_{2}=-\frac{1+\cos\varphi}{\cos\varphi}.
 \end{equation}
Using \eqref{f3}, \eqref{f8}, \eqref{ff6}, (\ref{ff7}) and (\ref{ff9})
 we find
 \begin{equation}
 A=\frac{1+3\cos\varphi}{\cos^{2}\varphi}\int^{\beta}_{\alpha} P(t,k_{1}t)dt.
 \end{equation}
  \item Similarly, if $F$ is on $c_{1}$ and $F$ acts on $c_{2}$ we get
  \begin{equation*}%\label{ff11}
 A=\frac{1+3\cos\varphi}{\cos^{2}\varphi}\int^{\beta}_{\alpha} P(t,k_{2}t)dt.
 \end{equation*}
\end{itemize}

\textbf{Case (D)} Finally, the condition $\varphi=\frac{\pi}{2}$ applied to \eqref{defw} yields $j=Qi$. Then $i$ and $j$ are isotropic vectors.
Therefore, from \eqref{f3} and \eqref{f8} it follows:
\begin{itemize}
\item $F=P(t,0)Qi$, $dr=(dt)i$. The work is  $A=\int^{\beta}_{\alpha}P(t,0)dt$.
\item $F=P(0,t)i$, $dr=(dt)Qi$. The work is
$A=\int^{\beta}_{\alpha}P(0,t)dt$.
\end{itemize}

\begin{table}
\caption{Work \textit{\textbf{A}} done by an isotropic vector force field $F$ along an isotropic curve}
\label{tab:1}       % Give a unique label
\begin{center}
\begin{tabular}{llll}
\noalign{\smallskip}\hline
$\varphi$ & $F$ acts on &  trajectory of $F$ & $A$ \\
\noalign{\smallskip}\hline
$(\arccos (-\frac{1}{3}), \frac{2\pi}{3})$ & - & no is. curves & - \\
\noalign{\smallskip}\hline
$\arccos(-\frac{1}{3})$ & $c:\ y=\sqrt{2}x$ & $c:\ y=\sqrt{2}x$ & 0\\
\noalign{\smallskip}\hline
$\big(0, \frac{\pi}{2}\big)\bigcup\big(\frac{\pi}{2}, \arccos(-\frac{1}{3})\big)$ & $c_{1}: y = k_{1}x$ & $c_{1}: y=k_{1}x$ & 0 \\
\noalign{\smallskip}\hline
$\big(0, \frac{\pi}{2}\big)\bigcup\big(\frac{\pi}{2}, \arccos(-\frac{1}{3})\big)$ & $c_{2}: y = k_{2}x$ & $c_{2}: y=k_{2}x$ & 0 \\
\noalign{\smallskip}\hline
$\big(0, \frac{\pi}{2}\big)\bigcup\big(\frac{\pi}{2}, \arccos(-\frac{1}{3})\big)$ & $c_{1}: y = k_{1}x$ & $c_{2}: y=k_{2}x$ & $A=\frac{1+3\cos\varphi}{\cos^{2}\varphi}\int^{\beta}_{\alpha} P(t,k_{1}t)dt$ \\
\noalign{\smallskip}\hline
$\big(0, \frac{\pi}{2}\big)\bigcup\big(\frac{\pi}{2}, \arccos(-\frac{1}{3})\big)$ & $c_{2}: y = k_{2}x$ & $c_{1}: y=k_{1}x$ & $A=\frac{1+3\cos\varphi}{\cos^{2}\varphi}\int^{\beta}_{\alpha} P(t,k_{2}t)dt$ \\
\noalign{\smallskip}\hline
$\frac{\pi}{2}$ & $c_{1}: x=0$ & $c_{2}: y=0$ &  $A=\int^{\beta}_{\alpha}P(t,0)dt$ \\
\noalign{\smallskip}\hline
$\frac{\pi}{2}$ & $c_{1}: y=0$ & $c_{2}: x=0$ &  $A=\int^{\beta}_{\alpha}P(0,t)dt$. \\
\noalign{\smallskip}\hline
\end{tabular}
\end{center}
\end{table}

The results in Case (A) -- Case (D) are summarized in Table \ref{tab:1}.

\section*{Acknowledgments}
This work is supported by project ''17-12 Supporting Intellectual Property'' of the Center for Research, Technology Transfer and Intellectual Property Protection, Agricultural University of Plovdiv, Bulgaria.

\end{document}